\documentclass[10pt]{amsart}

\usepackage{amssymb}

\newtheorem{thm1}{Theorem}[section]

\newtheorem{def1}[thm1]{Definition}

\newtheorem{cor1}[thm1]{Corollary}

\newtheorem{prop1}[thm1]{Proposition}
\newtheorem{ex1}[thm1]{Example}

\DeclareMathOperator{\supp}{supp} 
\DeclareMathOperator{\ind}{ind}

\begin{document}

\title[Minimal systems of generators]
{ Minimal systems of binomial generators and the indispensable
complex of a toric ideal}
\author[H. Charalambous]{Hara Charalambous}
\address { Department of Mathematics, Aristotle University of Thessaloniki, Thessaloniki 54124, GREECE
 } \email{hara@math.auth.gr}
\author[A. Katsabekis]{ Anargyros Katsabekis}
\address { Department of Mathematics, University of Ioannina,
Ioannina 45110, GREECE } \email{akatsabekis@in.gr}
\author[A. Thoma]{ Apostolos Thoma}
\address { Department of Mathematics, University of Ioannina,
Ioannina 45110, GREECE } \email{athoma@cc.uoi.gr}

\keywords{toric ideal, minimal systems of generators, indispensable
monomials, indispensable binomials} \subjclass{13F20, 05C99}

\begin{abstract}
\par Let $A=\{{\bf a}_1,\ldots,{\bf a}_m\} \subset \mathbb{Z}^n$ be a vector configuration and
$I_A \subset K[x_1,\ldots,x_m]$ its corresponding toric ideal. The
paper consists of two parts. In the first part we completely
determine the number of different minimal systems of binomial
generators of $I_A$. We also prove that generic toric ideals are
generated by indispensable binomials. In the second part we
associate to $A$ a simplicial complex $\Delta _{\ind(A)}$. We show
that the vertices of $\Delta _{\ind(A)}$ correspond to the
indispensable monomials of the toric ideal $I_A$, while one
dimensional facets of $\Delta _{\ind(A)}$ with minimal binomial
$A$-degree correspond to the indispensable binomials of $I_{A}$.
\end{abstract}
\thanks{This research was co-funded by the European Union
in the framework of the program ``Pythagoras" of the ``Operational
Program for Education and Initial Vocational Training" of the 3rd
Community Support Framework of the Hellenic Ministry of Education.
} \maketitle

\section{Introduction}

\par  Let $A=\{{\bf a}_1,\ldots,{\bf a}_m\}$
be a vector configuration in $\mathbb{Z}^n$ such that the affine
semigroup $\mathbb{N}A:=\{l_1{\bf a}_1+\cdots+l_m{\bf a}_m \ | \
l_i \in \mathbb{N}\}$ is pointed. Recall that $\mathbb{N}A$ is
{\em pointed} if zero is the only invertible element. Let $K$ be a
field of any characteristic; we grade the polynomial ring
$K[x_1,\ldots,x_m]$ by setting $\deg_{A}(x_i)={\bf a}_i$ for
$i=1,\ldots,m$. For ${\bf u}=(u_1,\ldots,u_m) \in \mathbb{N}^m$,
we define the $A$-{\em degree} of the monomial ${\bf x}^{{\bf
u}}:=x_1^{u_1} \cdots x_m^{u_m}$ to be \[ \deg_{A}({\bf x}^{{\bf
u}}):=u_1{\bf a}_1+\cdots+u_m{\bf a}_m \in \mathbb{N}A.\]  The
{\em toric ideal} $I_{A}$ associated to $A$ is the prime ideal
generated by all the binomials ${\bf x}^{{\bf u}}- {\bf x}^{{\bf
v}}$ such that $\deg_{A}({\bf x}^{{\bf u}})=\deg_{A}({\bf x}^{{\bf
v}})$ (see \cite{St}). For such binomials, we define $\deg_A({\bf
x}^{{\bf u}}- {\bf x}^{{\bf v}}):=\deg_{A}({\bf x}^{{\bf u}})$.
\par In general it is possible for a toric ideal $I_A$ to have more
than one minimal system of generators. We define $\nu (I_A)$ to be
the number of different minimal systems of binomial generators of
the toric ideal $I_A$, where the sign of a binomial does not
count. A recent problem arising from Algebraic Statistics, see
\cite{AT}, is when a toric ideal possesses a unique minimal system
of binomial generators, i.e. $\nu (I_A)=1$. To study this problem
Ohsugi and Hibi introduced in \cite{Hi-O} the notion of
indispensable binomials while Aoki, Takemura and Yoshida
introduced in \cite{ATY} the notion of indispensable monomials.
Moreover in \cite{Hi-O} a necessary and sufficient condition is
given for  toric ideals associated with certain finite graphs to
possess unique minimal systems of binomial generators. We recall
that a binomial $B={\bf x}^{{\bf u}}-{\bf x}^{{\bf v}} \in I_{A}$
is {\em indispensable} if every system of binomial generators of
$I_{A}$ contains $B$ or $-B$, while a monomial ${\bf x}^{{\bf u}}$
is {\em indispensable} if every system of binomial generators of
$I_{A}$ contains a binomial $B$ such that the ${\bf x}^{{\bf u}}$
is a monomial of $B$. \par In this article we use and extend
ideas-techniques developed by Diaconis, Sturmfels (see \cite{DS})
and Takemura, Aoki (see \cite{AT}) to study minimal systems of
generators of the toric ideal $I_A$ and also investigate the
notion of the indispensable complex of $I_A$, denoted by $\Delta
_{\ind(A)}$. In section \ref{numberofgeneratorssection}, we
construct graphs $G({\bf b})$, for every ${\bf b} \in
\mathbb{N}A$, and use them to provide a formula for $\nu (I_A)$.
We give criteria for a toric ideal to be generated by
indispensable binomials. In section 3 we determine a large class
of toric ideals, namely generic toric ideals, which have a unique
minimal system of binomial generators. In Section
\ref{indispcomplexsection} we define $\Delta _{\ind(A)}$ and we
show that this complex determines the indispensable monomials and
binomials. As an application we characterize principal toric
ideals in terms of $\Delta_{\ind(A)}$. \\

\section{The number of minimal generating sets of a toric ideal}\label{numberofgeneratorssection}

Let $A \subset \mathbb{Z}^n$ be a vector configuration so that
$\mathbb{N}A$ is pointed and $I_A \subset K[x_1,\ldots,x_m]$ its
corresponding toric ideal.  A vector ${\bf b}\in \mathbb{N}A$
 is called an {\em  Betti $A$-degree} if
$I_A$ has a minimal generating set containing an element of
$A$-degree ${\bf b}$. The Betti $A$-degrees are independent of the
choice of a minimal generating set of $I_A$, see \cite{CP, MS,
St}. The $A$-\emph{graded Betti number} $\beta_{0,{\bf b}}$ of
$I_A$ is the number of times $\bf b$ appears as  the $A$-degree of
a binomial in a given minimal generating set of $I_A$ and is also
an invariant of $I_A$.

The semigroup $\mathbb{N}A$ is pointed, so we can partially order
it with the relation
\[{\bf c} \geq {\bf d} \Longleftrightarrow \ \textrm{there is} \
{\bf e} \in \mathbb{N}A \ \textrm{such that} \ {\bf c}={\bf
d}+{\bf e}.\]
  For $I_A\neq \{0\}$ the minimal elements of the set
$\{\deg_{A}({\bf x}^{{\bf u}}) \ | \ {\bf x}^{{\bf u}}-{\bf
x}^{{\bf v}} \in I_{A}\} \subset \mathbb{N}A$ with respect to
$\geq$ are called {\em minimal binomial $A$-degrees}. Minimal
binomial $A$-degrees are always Betti $A$-degrees but the converse
is not true, as Example \ref{anexample} demonstrates. For any
${\bf b} \in \mathbb{N}A$ set \[I_{A,{\bf b}}:=({\bf x}^{{\bf
u}}-{\bf x}^{{\bf v}} \ | \ \deg_{A}({\bf x}^{{\bf
u}})=\deg_{A}({\bf x}^{{\bf v}}) \lneqq {\bf b}) \subset I_A.\]
\begin{def1} For a vector ${\bf b} \in \mathbb{N}A$ we define
$G({\bf b})$ to be the graph with vertices the elements of the
fiber
$$\deg_{A}^{-1}({\bf b})=\{{\bf x}^{{\bf u}} \ | \ \deg_{A}({\bf
x}^{{\bf u}})={\bf b}\}$$  and edges all the sets $\{{\bf x}^{{\bf
u}},{\bf x}^{{\bf v}}\}$ whenever ${\bf x}^{{\bf u}}-{\bf x}^{{\bf
v}} \in I_{A,{\bf b}}$.
\end{def1}
The fiber $\deg_{A}^{-1}({\bf b})$ has finitely many elements,
since the affine semigroup $\mathbb{N}A$ is pointed. If ${\bf
x}^{\bf u}$, ${\bf x}^{\bf v}$ are vertices of $G({\bf b})$ such
that $\gcd({\bf x}^{\bf u},{\bf x}^{\bf v}) \neq 1$, then $\{{\bf
x}^{\bf u},{\bf x}^{\bf v}\}$ is an edge of $G({\bf b})$. The next
proposition follows easily from the definition.

\begin{prop1}\label{connectedcomponents} Let ${\bf b} \in
\mathbb{N}A$. Every connected component of $G({\bf b})$ is a
complete subgraph.  The graph $G({\bf b})$ is not connected if and
only if ${\bf b}$ is a Betti $A$-degree.
\end{prop1}

\begin{ex1} \label{anexample} {\rm  Let \[
A=\{(2,2,2,0,0),(2,-2,-2,0,0),(2,2,-2,0,0),(2,-2,2,0,0),\]
\[ (3,0,0,3,3),(3,0,0,-3,-3),(3,0,0,3,-3)(3,0,0,-3,3)\}.\]
Using CoCoA, \cite{CoCoA} we see that  $I_A =(
x_1x_2-x_3x_4,x_5x_6-x_7x_8,x_1^3x_2^3-x_5^2x_6^2)$. The Betti
$A$-degrees are ${\bf b}_1=(4,0,0,0,0)$, ${\bf b}_2=(6,0,0,0,0)$
and ${\bf b}_3=(12,0,0,0,0)$. We note that ${\bf b}_3=2{\bf b}_2$,
so ${\bf b}_3$  is not a minimal binomial $A$-degree. The ideals
$I_{A,{\bf b}_1}$ and $I_{A,{\bf b}_2}$ are zero, while $I_{A,{\bf
b}_3}=(x_1x_2-x_3x_4,x_5x_6-x_7x_8)$. The graphs $G({\bf b})$ are
connected for all ${\bf b} \in \mathbb{N}A$ except for the Betti
$A$-degrees. In fact $G({\bf b}_1)$ and $G({\bf b}_2)$ consist of
two connected  components, $\{x_1x_2\}$ and $\{x_3x_4\}$ for
$G({\bf b}_1)$, $\{x_5x_6\}$ and $\{x_7x_8\}$ for $G({\bf b}_2)$,
while the connected components of $G({\bf b}_3)$ are
$\{x_1^3x_2^3, x_1^2x_2^2x_3x_4,\\ x_1x_2x_3^2x_4^2, x_3^3x_4^3\}$
and $\{x_5^2x_6^2, x_5x_6x_7x_8, x_7^2x_8^2\}$.}
\end{ex1}

Let $n_{{\bf b}}$ denote the number of  connected components of
$G({\bf b})$, this means that  \[ G({\bf
b})=\bigcup_{i=1}^{n_{{\bf b}}}G({\bf b})_{i} \] and $t_{i}({\bf
b})$ be the number of vertices of the $i$-component.  The next
proposition will be helpful in the sequel.

\begin{prop1} \label{prop} An $A$-degree ${\bf b}$
 is a minimal binomial $A$-degree if and only if every connected
component of $G({\bf b})$ is a singleton.
\end{prop1}
\noindent \textbf{Proof.} If ${\bf b}$ is a minimal binomial
$A$-degree, then $I_{A,{\bf b}}=\{0\}$ and every connected
component of $G({\bf b})$ is a singleton. Suppose now that ${\bf
b}$ is not minimal, i.e. ${\bf c} \lneqq {\bf b}$ for some minimal
binomial $A$-degree  $\bf c$. Thus there is a binomial $B={\bf
x}^{{\bf u}}-{\bf x}^{{\bf v}} \in I_A$, with $\deg_{A}(B)={\bf
c}$, and a monomial ${\bf x}^{{\bf a}}\neq 1$ such that ${\bf
b}={\bf c}+ \deg_{A}({\bf x}^{{\bf a}})$. Therefore ${\bf x}^{{\bf
a}+{\bf u}}$, ${\bf x}^{{\bf a}+{\bf v}}$ are vertices of $G({\bf
b})$ and belong to the same component of $G({\bf b})$ since ${\bf
x}^{{\bf a}+{\bf u}}-{\bf x}^{{\bf a}+{\bf v}}={\bf x}^{{\bf
a}}B\in I_{A,{\bf b}}$. \qed

Let $\mathcal{G} \subset I_A$ be a set of binomials. We recall the
definition of the graph $\Gamma({\bf b})_{\mathcal{G}}$,
\cite{DS}, and a  criterion  for $\mathcal{G}$ to be a generating
set of $I_A$, Theorem \ref{gentheorem}. Let $\Gamma({\bf
b})_{\mathcal{G}}$ be the graph with vertices the elements of
$\deg_{A}^{-1}({\bf b})$ and edges the sets $\{{\bf x}^{{\bf
u}},{\bf x}^{{\bf v}}\}$ whenever the binomial $$\frac{({\bf
x}^{{\bf u}}-{\bf x}^{{\bf v}})}{\gcd ( {\bf x}^{{\bf u}},{\bf
x}^{{\bf v}})} \ \ \textrm{or} \ \ \frac{({\bf x}^{{\bf v}}-{\bf
x}^{{\bf u}})}{\gcd ( {\bf x}^{{\bf u}},{\bf x}^{{\bf v}})}$$
belongs to $ \mathcal{G}$. In \cite{DS} the following theorem was
proved.

\begin{thm1}\label{gentheorem} \cite{DS} $\mathcal{G}$ is a generating set for
$I_A$ if and only if $\Gamma({\bf b})_{\mathcal{G}}$ is connected
for all ${\bf b}\in \mathbb{N}A$.
\end{thm1}
We consider the complete graph $\mathcal{S}_{{\bf b}}$ with
vertices the connected components $G({\bf b})_{ i}$ of $G({\bf
b})$, and we let $T_{\bf b}$ be a spanning tree of
$\mathcal{S}_{{\bf b}}$; for every edge of $T_{\bf b}$ joining the
components $G({\bf b})_i$ and $G({\bf b})_j$, we choose a binomial
${\bf x}^{\bf u}-{\bf x}^{\bf v}$ with ${\bf x}^{\bf u} \in G({\bf
b})_i$ and ${\bf x}^{\bf v} \in G({\bf b})_j$. We call
$\mathcal{F}_{T_{\bf b}}$ the collection of these binomials. Note
that if ${\bf b}$ is not a Betti $A$-degree, then
$\mathcal{F}_{T_{\bf b}}=\emptyset$.

\begin{thm1} \label{spanningtreetheorem}
The set $\mathcal{F}=\cup_{{\bf b} \in \mathbb{N}A}\mathcal{F}_{T_{\bf b}}$ is a
minimal generating set of $I_A$.
\end{thm1}
\noindent \textbf{Proof.} First we will prove that $\mathcal{F}$
is a generating set of $I_A$. From Theorem \ref{gentheorem} it is
enough to prove that $\Gamma({\bf b})_{\mathcal{F}}$ is connected
for every ${\bf b}$. We will prove the theorem by induction on
${\bf b}$. If $\bf b$ is a minimal binomial $A$-degree, the
vertices of $\Gamma({\bf b})_{\mathcal{F}}$, which are also the
vertices and the connected components of $G({\bf b})$, and the
tree $T_{\bf b}$ gives a path between any two vertices of $G({\bf
b})$. Next, let ${\bf b}$ be non-minimal binomial $A$-degree.
Suppose that $\Gamma({\bf b})_{\mathcal{F}}$ is connected for all
${\bf c} \lneqq {\bf b}$ and let ${\bf x}^{{\bf u}}$, ${\bf
x}^{{\bf v}}$ be two vertices of $\Gamma({\bf b})_{\mathcal{F}}$.
We will show that there is a path between these two vertices. We
will consider two cases, depending on whether the vertices are in
the same connected component of $G({\bf b})$ or not.

\begin{enumerate}
\item{} If ${\bf x}^{{\bf u}}$, ${\bf x}^{{\bf v}}$ are in the
same component $G({\bf b})_{i}$ of $G({\bf b})$, then ${\bf
x}^{{\bf u}}-{\bf x}^{{\bf v}}=\sum_{i} {\bf x}^{{\bf d}_i}({\bf
x}^{{\bf u}_i}-{\bf x}^{{\bf v}_i})$ where ${\bf x}^{{\bf u}_i}$,
${\bf x}^{{\bf v}_i}$ have $A$-degree ${\bf b}_i \lneqq {\bf b}$.
From the inductive hypothesis $\Gamma({\bf b}_i)_{\mathcal{F}}$ is
connected and there is a path from ${\bf x}^{{\bf u}_i}$ to ${\bf
x}^{{\bf v}_i}$. This gives a path from ${\bf x}^{{\bf d}_i}{\bf
x}^{{\bf u}_i}$ to ${\bf x}^{{\bf d}_i}{\bf x}^{{\bf v}_i}$ and
joining these paths we find a path from ${\bf x}^{{\bf u}}$ to
${\bf x}^{{\bf v}}$ in $\Gamma({\bf b})_{\mathcal{F}}$.

\item{} If  ${\bf x}^{{\bf u}}$, ${\bf x}^{{\bf v}}$ belong to
different components of $G({\bf b})$ we use the tree $T_{\bf b}$
to find a path between the two components. In each component we
use the previous case and/or the induction hypothesis to move
between vertices if needed. The join of these paths provides a
path from ${\bf x}^{{\bf u}}$ to ${\bf x}^{{\bf v}}$ in
$\Gamma({\bf b})_{\mathcal{F}}$.

\end{enumerate}

Next, we will show that no proper subset $\mathcal{F'}$ of
$\mathcal{F}$ generates $I_A$. Let $B={\bf x}^{{\bf u}}- {\bf
x}^{{\bf v}} \in \mathcal{F}\setminus \mathcal{F}'$, and
$\deg_{A}(B)={\bf b}$.  Since $B$ is an element of $
\mathcal{F}_{T_{\bf b}}$, it corresponds to an edge $\{ G({\bf
b})_{ i}, G({\bf b})_{j}\}$ of $T_{\bf b}$, and the monomials
${\bf x}^{{\bf u}}$, ${\bf x}^{{\bf v}}$ belong to different
components  of  $G({\bf b})$. Suppose that there was a path
$\{{\bf x}^{{\bf u}_1}= {\bf x}^{{\bf u}},{\bf x}^{{\bf
u}_2},\ldots,{\bf x}^{{\bf u}_t}={\bf x}^{{\bf v}}\}$ in
$\Gamma({\bf b})_{\mathcal{F'}}$ joining the vertices ${\bf
x}^{{\bf u}}$ and ${\bf x}^{{\bf v}}$. Certainly there are
monomials ${\bf x}^{{\bf u}_i}$, ${\bf x}^{{\bf u}_{i+1} }$ that
are in different connected components of $G({\bf b})$. Since
$\gcd({\bf x}^{{\bf u}_i},{\bf x}^{{\bf u}_{i+1}}) \neq 1$ implies
that the monomials ${\bf x}^{{\bf u}_i}$, ${\bf x}^{{\bf u}_{i+1}
}$ are in the same connected component of $G({\bf b})$, we
conclude that $\gcd({\bf x}^{{\bf u}_i},{\bf x}^{{\bf u}_{i+1}}) =
1$ for some $i$. In this case the binomial ${\bf x}^{{\bf
u}_i}-{\bf x}^{{\bf u}_{i+1} }$ is in $\mathcal{F'}$, has
$A$-degree ${\bf b}$, and it corresponds to an edge of $T_{\bf
b}$.  By considering these binomials and corresponding edges we
obtain a path in $T_{\bf b}$ joining the components $G({\bf b})_i,
G({\bf b})_j$ and not containing $\{G({\bf b})_i,G({\bf b})_j\}$
of $T_{\bf b}$. This is a contradiction since $T_{{\bf b}}$ is a
tree. \qed \\

The converse is also true; let $\mathcal{G}=\cup_{{\bf b} \in
\mathbb{N}A}\mathcal{G}_{{\bf b}}$ be a minimal generating set for
$I_A$ where $\mathcal{G}_{{\bf b}}$ consists of the binomials in
$\mathcal{G}$ of $A$-degree $\bf b$. We will show that
$\mathcal{G}_{{\bf b}}$ determines a spanning tree $T_{\bf b}$ of
$\mathcal{S}_{\bf b}$.

\begin{thm1}\label{generatorsspanningtree} Let $\mathcal{G}=\cup_{{\bf b} \in
\mathbb{N}A}\mathcal{G}_{{\bf b}}$ be a minimal generating set for
$I_A$. The binomials of $\mathcal{G}_{{\bf b}}$ determine a
spanning tree $T_{\bf b}$ of $\mathcal{S}_{\bf b}$.\end{thm1}

\noindent \textbf{Proof.} Let  $B={\bf x}^{{\bf u}}- {\bf x}^{{\bf
v}}\in \mathcal{G}_{{\bf b}}$. The monomials ${\bf x}^{{\bf u}}$,
${\bf x}^{{\bf v}}$ are in different connected components of
$G({\bf b})$, otherwise  $B$ is not a part of a minimal generating
set of $I_A$.  Therefore $B$ indicates an edge in
$\mathcal{S}_{\bf b}$. Let $T_{\bf b}$ be the union over $B\in
\mathcal{G}_{{\bf b}}$ of these edges. $T_{\bf b}$ is tree of
$\mathcal{S}_{\bf b}$, since  if $T_{\bf b}$ contains a cycle we
can delete a binomial from $\mathcal{G}$ and still generate the
ideal $I_A$, contradicting the minimality of $\mathcal{G}$.
Theorem \ref{gentheorem} guarantees that the tree $T_{\bf b}$ is
spanning. \qed \\

An immediate corollary of Theorems \ref{spanningtreetheorem} and
\ref{generatorsspanningtree} concerns the indispensable monomials.

\begin{cor1} ${\bf x}^{\bf u}$ is an indispensable monomial of A-degree ${\bf
b}$ if and only if $\{ {\bf x}^{\bf u} \}$ is a component of
$G({\bf b})$. \end{cor1}

We use Theorems \ref{spanningtreetheorem} and
\ref{generatorsspanningtree} to compute $\nu (I_A)$, the number of
minimal generating sets of $I_A$. For each ${\bf b} \in
\mathbb{N}A$ the number of possible spanning trees $T_{\bf b}$
depends on $n_{{\bf b}}$, the number of connected components of
$G({\bf b})$. For a given spanning tree $T_{\bf b}$ the number of
possible binomial sets $\mathcal{F}_{T_{\bf b}}$ (up to a sign)
depends on $t_{i}({\bf b})$, the number of vertices of  $G({\bf
b})_{ i}$. These numbers determine $\nu (I_A)$. We note that the
sum $t_{1}({\bf b})+\cdots +t_{n_{\bf b}}({\bf b} )$ is equal to
$| \deg_{A}^{-1}({\bf b})|$, the cardinality of the fiber set
$\deg_{A}^{-1}({\bf b})$. We also point out that
$|\mathcal{F}_{T_{\bf b}}|= n_{\bf b}-1$ and that
$|\mathcal{F}_{T_{\bf b}}|=\beta _{0,{\bf b}}$, the $A$-graded
Betti number of $I_A$.

\begin{thm1} \label{number} For a toric ideal $I_A$ we have that
\[ \nu (I_A)=\prod_{{\bf b}\in \mathbb{N}A} t_{1}({\bf b})\cdots t_{n_{\bf b}}({\bf
b})(t_{1}({\bf b})+\cdots +t_{n_{\bf b}}({\bf b} ))^{n_{\bf b}-2}\]
where $n_{{\bf b}}$ is the number of connected components of $G({\bf
b})$ and $t_{i}({\bf b})$ is the number of vertices of the connected
component $G({\bf b})_{i}$ of the graph $G({\bf b})$.
\end{thm1}

\noindent \textbf{Proof.} Let $d_{i}$ be the degree of $G({\bf
b})_{i}$ in a spanning tree $T_{\bf b}$, i.e. the number of edges
of $T_{\bf b}$ incident with $G({\bf b})_{i}$. We have that $\sum
_{i=1}^{n_{\bf b}}d_i=2n_{\bf b}-2$. There are
\[{\frac{(n_{\bf b}-2)!}{(d_1-1)!(d_2-1)!\cdots
(d_{n_{\bf b}}-1)!}}\] such spanning trees, see for example the
proof of Cayley's formula in \cite{LW}. For fixed $T_{\bf b}$ with
degrees $d_i$, there are $(t_{i}({\bf b}))^{d_i}$ choices for the
monomials for the edges involving the vertex $G({\bf b})_{i}$.
This implies that the number of possible binomial sets
$\mathcal{F}_{T_{\bf b}}$ is $(t_{1}({\bf b}))^{d_1}\cdots
(t_{n_{\bf b}}({\bf b}))^{d_{n_{\bf b}}}$. Therefore the total
number of all possible $\mathcal{F}_{T_{\bf b}}$ is
\[ \sum_{d_1+\dots +d_{n_{\bf b}}=2n_{\bf b}-2}{\frac{(n_{\bf
b}-2)!}{(d_1-1)!(d_2-1)!\cdots (d_{n_{\bf b}}-1)!}}(t_{1}({\bf
b}))^{d_1}\cdots (t_{n_{\bf b}}({\bf b}))^{d_{n_{\bf b}}}=\]
\[=t_{1}({\bf b})\cdots t_{n_{\bf b}}({\bf b})(t_{1}({\bf
b})+\cdots +t_{n_{\bf b}}({\bf b}))^{n_{\bf b}-2}.\]
 \qed

We point out that if $t_i({\bf b})=1$ for all $i$, then the number
of possible spanning trees is ${n_{{\bf b}}}^{n_{{\bf b}}-2}$, (
Cayley's formula, see \cite{Ca}). We also note that if $n_{\bf
b}=1$, for some ${\bf b}\in \mathbb{N}A$,  then the factor
$t_{1}({\bf b})(t_{1}({\bf b}))^{-1}$ in the above product has
value 1. Thus the contributions to $\nu (I_A)$ come only from
Betti $A$-degrees ${\bf b} \in \mathbb{N}A$. On the other hand we
have a unique choice for a generator of degree $\bf b$ when
$n_{\bf b}= 2$ and $t_1({\bf b})=t_2({\bf b})=1$. Thus in these
cases $G({\bf b})$ consists of two isolated vertices and by
Proposition \ref{prop}, $\bf b$ is minimal. These remarks prove
the following:

\begin{cor1} \label{corind} Let $B={\bf
x}^{{\bf u}}-{\bf x}^{{\bf v}} \in I_A$ with $A$-degree ${\bf b}$.
$B$ is indispensable if and only if the graph $G({\bf b})$
consists of two connected components, $\{{\bf x}^{{\bf u}}\}$ and
$\{{\bf x}^{{\bf v}}\}$. Moreover  ${\bf b}$ is minimal binomial
$A$-degree.
\end{cor1}

\begin{cor1} \label{numgencor} Suppose that the Betti $A$-degrees ${\bf b}_1,\ldots,{\bf b}_q$ of
 $I_A$ are minimal binomial
$A$-degrees. Then
\[\nu
(I_A)=(\beta_{0,{\bf b}_1}+1)^{\beta_{0,{\bf b}_1}-1} \cdots
(\beta_{0,{\bf b}_q}+1)^{\beta_{0,{\bf b}_q}-1}.\]
\end{cor1}

\noindent \textbf{Proof.} By Proposition \ref{prop}, the connected
components of $G({\bf b}_i)$ are singletons. It follows that
$t_j({\bf b}_i)=1$ and that $n_{{\bf b}_i}= \sum t_j({\bf b}_i)$.
Moreover $\beta _{0,{\bf
b}_i}=|\mathcal{F}_{{\bf b}_i}|=n_{{\bf b}_i}-1$. \qed\\

The next theorem provides a necessary and sufficient condition for
a toric ideal to be generated by its indispensable binomials. It
is a generalization of Corollary 2.1 in \cite{AT}.

\begin{thm1} \label{indispenscor} The ideal $I_A$ is generated by its indispensable binomials if
and only if the Betti $A$-degrees ${\bf b}_1,\ldots,{\bf b}_q$  of
$I_A$ are minimal binomial $A$-degrees  and $\beta _{0,{\bf
b}_i}=1$.
\end{thm1}

\noindent \textbf{Proof.} Suppose that $I_A$ is generated by
indispensable binomials, then $\nu (I_A)=1$ and therefore, from
Theorem \ref{number}, $t_j({\bf b}_i)=1$ and $n_{{\bf b}_i}=2$,
for all $j, i$. Thus $\beta _{0,{\bf b}_i}=1$. Now Proposition
\ref{prop} together with the fact that $t_j({\bf b}_i)=1$ implies
that all ${\bf b}_i$ are minimal binomial A-degrees. \qed\\

We point out that the above theorem implies that in the case that
a toric ideal $I_A$ is generated by indispensable binomials no two
minimal generators can have the same $A$-degree. We compute $\nu
(I_A)$ in the following example.

\begin{ex1} \label{basicexample}{\rm Let
$A=\{a_0=k,a_1=1,\ldots,a_n=1\} \subset \mathbb{N}$ be a set of
$n+1$ natural numbers with $k>1$ and  $I_A \subset
K[x_0,x_1,\ldots,x_n]$, the corresponding toric ideal. The ideal
$I_A$ is minimally generated by the binomials $x_0-x_1^k, x_1-x_2,
\ldots, x_{n-1}-x_n$. The Betti $A$-degrees are ${\bf b}_1=1$ and
${\bf b}_2=k$, while the $A$-graded Betti numbers are
$\beta_{0,1}=n-1$ and $\beta_{0,k}=1$. Also $G(1)$ consists of $n$
vertices, each one being a connected component, and $G(k)$ has two
connected components, the singleton $\{ x_0 \}$ and the complete
graph on the ${k+n-1}\choose {n-1}$ vertices $x_1^k,
x_1^{k-1}x_2,\ldots, x_n^k$. Thus
\[\nu
(I_A)=n^{n-2}{ {k+n-1}\choose {n-1}}. \]}\\
\end{ex1}

\section{Generic toric ideals are generated by indispensable binomials}

Generic toric ideals were introduced in \cite{PS} by Peeva and
Sturmfels. The term generic is justified due to a result of Barany
and Scarf (see \cite{BS}) in integer programming theory which
shows that, in a well defined sense, almost all toric ideals are
generic. Given a vector ${\bf
\alpha}=(\alpha_{1},\ldots,\alpha_{m}) \in \mathbb{N}^m$, the {\em
support} of ${\bf \alpha}$, denoted by $\supp({\bf \alpha})$, is
the set $\{i \in \{1,\ldots,m\} \ | \ \alpha_i \neq 0\}$. For a
monomial ${\bf x}^{{\bf u}}$ we define $\supp({\bf x}^{{\bf
u}}):=\supp({\bf u})$. A toric ideal $I_A \subset
K[x_1,\ldots,x_m]$ is called { \em generic} if it is minimally
generated by binomials with full support, i.e.,
\[ I_A=({\bf
x}^{{\bf u}_1}- {\bf x}^{{\bf v}_1},\ldots,{\bf x}^{{\bf u}_r}-
{\bf x}^{{\bf v}_r})\] where $\supp({\bf u}_i) \cup \supp({\bf
v}_i)=\{1,\ldots,m\}$ for every $i \in \{1,\ldots,r\}$, see
\cite{PS}. We will  prove that the minimal binomial generating set
of $I_A$ is a unique.

\begin{thm1}\label{genericthm} If $I_A$ is a generic toric ideal, then $\nu
(I_A)=1 $ and  $I_A$ is generated by its indispensable binomials.
\end{thm1}

\noindent \textbf{Proof.} Let $\{B_1, B_2,\ldots ,B_s\}$ be a
minimal generating set of $I_A$ of full support where $B_i= {\bf
x}^{{\bf u}_i}-{\bf x}^{{\bf v}_i}$ and ${\bf b}_i=
\deg_{A}(B_i)$. We will show that all ${\bf b}_i$ are minimal
binomial $A$-degrees. Suppose that one of them, say ${\bf b}_1$ is
not minimal and that  ${\bf b}_j+{\bf b}={\bf b}_1$ for ${\bf
b}\in \mathbb{N}A$. It follows that ${\bf x}^{{\bf a}}{\bf
x}^{{\bf u}_j}$, ${\bf x}^{{\bf a}}{\bf x}^{{\bf v}_j}$ are in the
same connected component of $G({\bf b}_1)$, where ${\bf x}^{{\bf
a}}$ is  a monomial of A-degree ${\bf b}$. Since $\supp({\bf u}_j)
\cup \supp({\bf v}_j)=\{1,\ldots,m\}$ it follows that at least one
of $\gcd({\bf x}^{{\bf u}_1},{\bf x}^{{\bf a}}{\bf x}^{{\bf u}_j})
$ or $\gcd({\bf x}^{{\bf u}_1},{\bf x}^{{\bf a}}{\bf x}^{{\bf
v}_j})$ is not $1$ and ${\bf x}^{{\bf u}_1}$ belongs to the same
connected component of $G({\bf b}_1)$ as ${\bf x}^{{\bf a}}{\bf
x}^{{\bf u}_j}$ and ${\bf x}^{{\bf a}}{\bf x}^{{\bf v}_j}$. The
same holds for ${\bf x}^{{\bf v}_1}$. This is a contradiction
since ${\bf x}^{{\bf u}_1}$ and ${\bf x}^{{\bf v}_1}$ belong to
different connected components of $G({\bf b}_1)$.

Next we will show that  $\beta_{0,{\bf b}_i}=1$. Suppose that one
of them, say  $\beta_{0,{\bf b}_1}=|\mathcal{F}_{T_{{\bf b}_1}}|$
is greater than 1. Since  ${\bf b}_1$ is minimal, the connected
components of $G({\bf b}_1)$ are singletons and $n({\bf b}_1)\ge
3$. It follows that in ${T_{{\bf b}_1}}$, two edges share a
vertex, and in $\mathcal{F}_{T_{{\bf b}_1}}$ the same monomial
will necessarily appear in two of our binomial generators. Since
the generators have full support the other two monomials have the
same support, thus they have a nontrivial common factor and are in
the same component of $G({\bf b}_1)$, a contradiction. \qed

\begin{ex1} {\rm Consider the vector configuration
$A=\{20,24,25,31\}$. Using CoCoA, \cite{CoCoA}, we see that
$I_A=(x_3^3-x_1x_2x_4,
x_1^4-x_2x_3x_4,x_4^3-x_1x_2^2x_3,x_2^4-x_1^2x_3x_4,x_1^3x_3^2-x_2^2x_4^2,x_1^2x_2^3-x_3^2x_4^2,
x_1^3x_4^2-x_2^3x_3^2)$ and therefore it is a generic ideal. We note
that the variety $V(I_A)$ is the generic monomial curve in the
affine $4$-space $A^4$ of smallest degree, see Example 4.5 in
\cite{PS}. By Theorem \ref{genericthm} $v(I_A)=1$ and the above
generators are indispensable binomials. It follows that the minimal
binomial $A$-degrees are
$75$, $80$, $93$, $96$, $110$, $112$ and $122$. }\\
\end{ex1}

\section{The indispensable complex of a vector configuration} \label{indispcomplexsection}

Consider a vector configuration $A=\{{\bf a}_1,\ldots,{\bf a}_m\}$
in $\mathbb{Z}^n$ with $\mathbb{N}A$ pointed, and the toric ideal
$I_{A} \subset K[x_1,\ldots,x_m]$. In \cite{H-O} it is proved that
a binomial $B$ is indispensable if and only if either $B$ or $-B$
belongs to the reduced Gr\"{o}bner base of $I_A$ for any
lexicographic term order on $K[x_1,\ldots,x_m]$. In \cite{ATY} it
is shown that a monomial $M$ is indispensable if the reduced
Gr\"{o}bner base of $I_A$, with respect to any lexicographic term
order on $K[x_1,\ldots,x_m]$, contains a binomial $B$ such that
 $M$ is a monomial of $B$. We are going to provide a more
efficient way to check if a binomial is indispensable and
respectively for a monomial. Namely we will give a criterion that
provides the indispensable binomials and monomials with only the
information from one specific generating set of $I_A$.

We let $\mathcal{M}_{A}$ be the monomial ideal generated by all
 ${\bf x}^{{\bf u}}$ for which there exists a nonzero
 ${\bf x}^{{\bf u}}-{\bf x}^{{\bf v}} \in I_{A}$; in other
words given a vector ${\bf u}=(u_1,\ldots,u_m) \in \mathbb{N}^m$
the monomial ${\bf x}^{{\bf u}}$ belongs to $\mathcal{M}_{A}$ if
and only if there exists ${\bf v}=(v_1,\ldots,v_m) \in
\mathbb{N}^m$ such that ${\bf v} \neq {\bf u}$, i.e. $v_i \neq
u_i$ for some $i$, and $\deg_A ({\bf x}^{\bf u})=\deg_A ({\bf
x}^{\bf v})$. We note  that if $\{B_1={\bf x}^{{\bf u}_{1}}- {\bf
x}^{{\bf v}_{1}},\ldots, B_s= {\bf x}^{{\bf u}_{s}}- {\bf x}^{{\bf
v}_{s}}\}$ is a generating set of $I_A$ then $\mathcal{M}_{A}=
({\bf x}^{{\bf u}_{1}},\ldots, {\bf x}^{{\bf u}_{s}}, {\bf
x}^{{\bf v}_{1}},\ldots, {\bf x}^{{\bf v}_{s}})$. Let
$T_A:=\{M_1,\ldots,M_k\}$ be the unique minimal monomial
generating set of $\mathcal{M}_{A}$.

\begin{prop1}\label{indmon}  The indispensable monomials of $I_A$ are precisely
the elements of $T_A$.
\end{prop1}

\noindent \textbf{Proof.} First we will prove that the elements of
$T_A$ are indispensable monomials. Let $\{B_1,\ldots,B_s\}$ be a
minimal generating set of $I_{A}$. Set $M_j:={\bf x}^{{\bf u}}$
for $j \in \{1,\ldots,k\}$. Since ${\bf x}^{{\bf u}}-{\bf x}^{{\bf
v}}$ is in $I_{A}$ for some $\bf v$, it follows that there is an
$i\in \{1,\ldots,s\}$ and a monomial $N$ of $B_i$ such that $N$
divides ${\bf x}^{{\bf u}}$ and thus ${\bf x}^{{\bf u}}=N$.

Conversely consider an indispensable monomial ${\bf x}^{{\bf u}}$
of $I_A$ and assume that is not an element of $T_A$ then  ${\bf
x}^{{\bf u}}=M_j {\bf x}^{{\bf c}}$ for some $j\in \{1,\ldots,
k\}$ and ${\bf c} \neq {\bf 0}$. By our previous argument $M_j$ is
indispensable. Without loss of generality we may assume that
$B_1=M_j-{\bf x}^{{\bf z}}$. If $B_j={\bf x}^{{\bf u}}-{\bf
x}^{{\bf v}}$, then
\[ B_j':={\bf x}^{{\bf c}}{\bf x}^{{\bf z}}-{\bf x}^{{\bf
v}}= B_j-{\bf x}^{{\bf c}} B_1 \in I_A \] and therefore
$I_A=(B_1,\ldots,B_{j-1}, B_j',B_{j+1}, \ldots,B_s)$. This way we
can eliminate ${\bf x}^{{\bf u}}$ from all the elements of the
generating set of $I_A$, a contradiction to the fact that ${\bf
x}^{{\bf u}}$ is indispensable. \qed

\begin{def1} {\rm We define the indispensable complex  $\Delta _{\ind(A)}$ to be
the simplicial complex with vertices the elements of $T_A$ and
faces all subsets of $T_A$ consisting of monomials with  the same
$A$-degree.}
\end{def1}
By Proposition \ref{indmon} the indispensable monomials are the
vertices of $\Delta _{\ind(A)}$. The connected components consist
of the vertices of the same $A$-degree and are simplices of
$\Delta _{\ind(A)}$, actually are facets. Different connected
components have different $A$-degrees. We compute $\Delta
_{\ind(A)}$ in the following example.

\begin{ex1} \rm{ In Example \ref{basicexample} we have that $\mathcal{M}_{A}=(x_0, x_1,\ldots, x_n)$
and also the facets of $\Delta _{\ind(A)}$ are $\{x_0\}$ and
$\{x_1, \ldots, x_n\}$. }
\end{ex1}

It follows easily that whenever $\deg_{A}({\bf x}^{{\bf u}})$ is a
minimal binomial $A$-degree, then ${\bf x}^{{\bf u}}\in T_A$. The
converse is not true in general. Indeed in Example
\ref{basicexample},  $x_0$ belongs to $T_A$ but  $\deg_{A}(x_0)$
is not minimal. Next we give a criterion that determines the
indispensable binomials.

\begin{thm1}\label{indispensablebinom}  A binomial $B={\bf x}^{{\bf u}}-{\bf
x}^{{\bf v}} \in I_{A}$ is indispensable if and only if $\{{\bf
x}^{{\bf u}}, {\bf x}^{{\bf v}}\}$ is a 1-dimensional facet of
$\Delta _{\ind(A)}$ and $\deg_{A}(B)$ is a minimal binomial
$A$-degree. \end{thm1}

\noindent \textbf{Proof.} Let ${\bf b}=\deg_{A}(B)$. Suppose that
$\{{\bf x}^{{\bf u}}, {\bf x}^{{\bf v}}\}$ is a 1-dimensional
facet of $\Delta _{\ind(A)}$ and  ${\bf b}$ is minimal binomial
$A$-degree. By Proposition \ref{prop}, minimality of $\bf b$
implies that the elements of $\deg_A^{-1}({\bf b})$ which are the
vertices of $G({\bf b})$, are vertices of $\Delta _{\ind(A)}$ and
the connected components of $G({\bf b})$ are singletons. Since
$\Delta _{\ind(A)}$ contains only two vertices of $A$-degree ${\bf
b}$, $G({\bf b})$ consists of two connected components, $\{{\bf
x}^{{\bf u}}\}$ and $\{{\bf x}^{{\bf v}}\}$ and $B$ is
indispensable by  Corollary \ref{corind}. The other direction is
done by reversing the implications. \qed

Theorem \ref{indispensablebinom} shows that the toric ideal $I_A$
of Example \ref{basicexample} has no indispensable binomials for
$n>2$. Indeed in this case the indispensable complex of $I_A$
contains no $1$-simplices that are facets.

We remark that to check the minimality of the $A$-degree $\bf b$
of the binomial $B\in I_A$ it is enough to compare $\bf b$ with
the $A$-degrees of the vertices of $\Delta _{\ind(A)}$. Thus given
any generating set of $I_A$ one can compute $T_A$ and construct
the simplicial complex $\Delta _{\ind(A)}$. The elements of $T_A$
are the indispensable monomials and the 1-dimensional facets of
$\Delta _{\ind(A)}$ of minimal binomial A-degree are the
indispensable binomials.

\begin{ex1} \label{exampindinspbin}{\rm Let
\[ A=\{(2,1,0),(1,2,0),(2,0,1),(1,0,2),(0,2,1),(0,1,2)\}.\]
Using CoCoA, \cite{CoCoA},  we see that $I_A$ is minimally
generated by: $x_1x_6-x_2x_4, x_1x_6-x_3x_5, x_2^2x_3 - x_1^2x_5,
x_2x_3^2 - x_1^2x_4, x_1x_5^2 - x_2^2x_6, x_1x_4^2 - x_3^2x_6,
x_4^2x_5 - x_3x_6^2, x_1x_4x_5 - x_2x_3x_6, x_4x_5^2 - x_2x_6^2$.
Moreover $T_A=\{ M_1=x_1x_6,\ M_2=x_2x_4,\ M_3=x_3x_5,\
M_4=x_3x_6^2,\ M_5=x_4^2x_5,\ M_6=x_2x_6^2,\  M_7=x_4x_5^2,\
M_8=x_3^2x_6,\ M_9=x_1x_4^2, \ M_{10}=x_2^2x_6,\ M_{11}=x_1x_5^2,\
M_{12}=x_1^2x_5,\ M_{13}=x_2^2x_3,\ M_{14}=x_1^2x_4,\
M_{15}=x_2x_3^2, M_{16}=x_2x_3x_6,\ M_{17}=x_1x_4x_5 \}$. It
follows that $\Delta _{\ind(A)}$ is a simplicial complex on $17$
vertices and its connected components are the facets
\[ \{M_1, M_2, M_3\}, \{M_4, M_5\}, \{M_6, M_7\}, \{M_8, M_9\}, \]
\[\{M_{10}, M_{11}\}, \{M_{12}, M_{13}\}, \{M_{14}, M_{15}\}, \{M_{16}, M_{17}\}.\]
The  $A$-degrees of the components are accordingly
\[ (2, 2, 2),  (2, 2, 5),  (1, 4, 4),  (4, 1, 4),\] \[  (2, 5, 2),  (4, 4, 1),  (5, 2, 2),
(3, 3, 3).\] All of them are minimal binomial $A$-degrees and thus
$I_A$ has 7 indispensable binomials corresponding to the
$1$-dimensional facets. We see that all non zero $A$-graded Betti
numbers equal to $1$, except from $\beta_{0,(2,2,2)}$ which equals
to $2$. From Corollary \ref{numgencor} we take that $\nu (I_A)=3$.}
\end{ex1}

The next corollary gives a necessary condition for a toric ideal
to be generated by the indispensable binomials.
\begin{cor1} \label{componentsindispensable} Let $A=\{{\bf a}_1,\ldots,{\bf a}_m\}$ be a vector configuration in
$\mathbb{Z}^n$. If $I_A$ is generated by the indispensable
binomials, then every connected component of $\Delta _{\ind(A)}$
is $1$-simplex.
\end{cor1}

\noindent \textbf{Proof.} Let $\{ B_1,\ldots,B_s\} $ be  a minimal
generating set of $I_A$ consisting of indispensable binomials
$B_i={\bf x}^{{\bf u}_i}-{\bf x}^{{\bf v}_i} $. We note that the
monomials of the $B_i$ are all indispensable and form $T_A$. Thus
if a face of $\Delta _{\ind(A)}$ contains $ {\bf x}^{{\bf u}_i} $
it also contains $ {\bf x}^{{\bf v}_i}$. By Theorem
\ref{indispensablebinom} $\{{\bf x}^{{\bf u}},{\bf x}^{{\bf v}}\}$
is a facet of $\Delta _{\ind(A)}$. \qed\\

The next example shows that the converse of Corollary
\ref{componentsindispensable} does not hold.

\begin{ex1} {\rm  We return to Example \ref{anexample}. The
simplicial complex $\Delta _{\ind(A)}$ consists of only  two
1-simplices $\{x_1x_2,x_3x_4\}$, $\{x_5x_6,x_7x_8\}$, the
indispensable binomials are $x_1x_2-x_3x_4$, $x_5x_6-x_7x_8$ and
$I_A \neq (x_1x_2-x_3x_4,x_5x_6-x_7x_8)$.}
\end{ex1}

When $\Delta _{\ind(A)}$ is a $1$-simplex the next proposition
shows that $I_A$ is principal and therefore generated by an
indispensable binomial.

\begin{prop1} The simplicial complex $\Delta _{\ind(A)}$ is a $1$-simplex if and only if $I_A$
is a principal ideal.
\end{prop1}

\noindent \textbf{Proof.} One direction of this Proposition is
trivial. For the converse assume that  $\Delta _{\ind(A)}=\{{\bf
x}^{{\bf u}_1},{\bf x}^{{\bf v}_1}\}$ and let $B_1:={\bf x}^{{\bf
u}_1}-{\bf x}^{{\bf v}_1}$. We will show that $I_A=(B_1)$. Let $B=
{\bf x}^{{\bf u}}-{\bf x}^{{\bf v}}$ be the binomial of minimal
binomial $A$-degree such that $B\in I_A \setminus (B_1)$. Since
${\bf x}^{{\bf u}}={\bf x}^{{\bf c}}{\bf x}^{{\bf u}_1}$ and ${\bf
x}^{{\bf v}}={\bf x}^{{\bf d}}{\bf x}^{{\bf v}_1}$, where ${\bf
x}^{{\bf c}}\neq {\bf x}^{{\bf d}}$, and none of them equals to 1,
we have that

\[ {\bf x}^{{\bf v}_1}({\bf x}^{{\bf c}}-{\bf x}^{{\bf d}})= {\bf x}^{{\bf
c}}B_1-B.\] Therefore $0\neq {\bf x}^{{\bf c}}-{\bf x}^{{\bf d}}
\in I_A$, while $\deg_A ({\bf x}^{{\bf c}}) \lneqq \deg_{A}({\bf
x}^{{\bf
u}})$ a contradiction. \qed\\

\end{document}